\documentclass[leqno,12pt]{amsart} %leqno is the option to put formula numbers on the left side
\usepackage[top=30truemm,bottom=30truemm,left=25truemm,right=25truemm]{geometry}
\usepackage{amssymb}
\usepackage{amsmath}
\usepackage{amsthm}
\usepackage{amscd}
\usepackage{mathrsfs}
\usepackage{graphicx}
\usepackage[dvips]{color}
\usepackage[all]{xy}
\usepackage{url}
  \makeatletter

\@addtoreset{equation}{section}
\makeatother
\theoremstyle{plain} %text of this environment is typesetted in italics
\newtheorem{theorem}{\indent\bf Theorem}[section]
\newtheorem{lemma}[theorem]{\indent\bf Lemma}

\newtheorem{proposition}[theorem]{\indent\bf Proposition}

\newtheorem{conjecture}[theorem]{\indent\bf Conjecture}
\theoremstyle{definition} %text of this environment is typesetted in roman letters
\newtheorem{definition}[theorem]{\indent\bf Definition}

\newtheorem{problem}[theorem]{\indent\bf Problem}

\newcommand{\ddbar}{\sqrt{-1} \partial \overline{\partial}}
\newcommand{\dbar}{\overline{\partial}}
\newcommand{\ai}{\sqrt{-1}}
\newcommand{\phinyoro}{\widetilde{\varphi}}
\newcommand{\llangle}{\langle\!\langle}
\newcommand{\rrangle}{\rangle\!\rangle}
\newcommand{\tr}{{\rm tr}}

\begin{document}
\pagestyle{plain}
\thispagestyle{plain}

\title[]
{From H\"ormander's $L^2$-estimates to partial positivity}

\author[T. INAYAMA]{Takahiro INAYAMA}
\address{Graduate School of Mathematical Sciences\\
The University of Tokyo\\
3-8-1 Komaba, Meguro-ku\\
Tokyo, 153-8914\\
Japan
}
\email{inayama@ms.u-tokyo.ac.jp}
\email{inayama570@gmail.com}
%\date{\today}

%%%%%%%%%%%%%%%%%%%%%%%%%%%%%%%%%%%%%%%%%%%%%%%%%%%%%%
\begin{abstract}
%In this article, we propose a new approach to 
In this article, using a twisted version of H\"ormander's $L^2$-estimate, we give new characterizations of notions of partial positivity, which are uniform $q$-positivity and RC-positivity.
%As an application of these results, we study a generalization of  the Pr\'ekopa-Berndtsson theorem. 
%Moreover, we show the RC-positivity of some sort of direct image bundle.
We also discuss the definition of uniform $q$-positivity for singular Hermitian metrics.
\end{abstract}

%%%%%%%%%%%%%%%%%%%%%%%%%%%%%%%%%%%%%%%%%%%%%%%%%%%%%%

\maketitle
\setcounter{tocdepth}{2}
%\tableofcontents
%%%%%%%%%%%%%%%%%%%%%%%%%%%%%%%%%%%%%%%%%%%%%%%%%%%%%%

\section{Introduction}

In this article, we give a new characterization of partial positivity, which is called uniform $q$-positivity (cf. Definition \ref{def:qpositive}) via H\"ormander's $L^2$-estimate.
The statement is the following. 

\begin{theorem}\label{thm:qpositive}
	Let $D$ be a bounded pseudoconvex domain in $\mathbb{C}^n_{ z}$, $\omega=\ddbar |z|^2$ be the standard K\"ahler metric on $D$, and $L\to D$ be a line bundle over $D$. For a smooth Hermitian metric $h$ on $L$ and a non-negative constant $c\geq 0$ on $D$, the following properties are equivalent for $1\leq q\leq n$: 
	\begin{enumerate}
		\item The summation of any distinct $q$ eigenvalues $($counting multiplicity$)$ of the Chern curvature $\ai \Theta_{(L, h)}$ of $(L, h)$ with respect to $\omega$ is greater than or equal to $c$. 
		\item For any smooth strictly plurisubharmonic function $\psi$ and any smooth $\dbar$-closed $L$-valued $(n,q)$-form $f$ with compact support, 
		%$$
		%\int_{D} \langle [\ddbar \psi, \Lambda_\omega]+c)^{-1}
		%f, f \rangle_{(\omega, h)} e^{-\psi}dV_\omega < +\infty,
		%$$
		there exists $L$-valued $(n,q-1)$-form $u$ satisfying 
		$\dbar u=f$ and 
		$$
		\int_{D}|u|^2_{(\omega, h)}e^{-\psi} dV_\omega \leq \int_{D} \langle ([\ddbar \psi, \Lambda_\omega]+c)^{-1}f, f\rangle_{(\omega, h)} e^{-\psi}dV_\omega.
		$$ 
	\end{enumerate}
\end{theorem}

Condition (2) in Theorem \ref{thm:qpositive} allows us to add a weight $\psi$. 
Taking an arbitrary weight, we can estimate the curvature $\ai\Theta_{(L, h)}$. 
This type of condition was firstly introduced in \cite{HI20}, where it was called the twisted H\"ormander condition. 
%In the paper, the authors proved that this condition implies Griffiths semi-positivity under some regularity conditions.
Next, in \cite{DNW19} and \cite{DNWZ20}, Deng et al. generalized this notion and introduced the optimal $L^p$-estimate condition, which corresponded to the particular case of the twisted H\"ormander condition when $p=2$.
These studies provide new characterizations of positivity based on the H\"ormander-type condition, which was initially observed by Berndtsson in \cite{Ber98}. 
%
%
%Roughly speaking, this work can be regarded as a converse of H\"ormander's $L^2$-estimate, which was initially observed by Berndtsson in \cite{Ber98}. 
%After that, in \cite{DNW19} and \cite{DNWZ20}, the authors introduced the optimal $L^p$-estimate condition, which corresponded to the particular case of the twisted H\"ormander condition when $p=2$.
%Moreover, by using the Bochner-Kodaira-Nakano identity and developing an excellent method, Deng et al. proved that the optimal $L^2$-estimate condition implies (Nakano) semi-positivity for smooth Hermitian metrics on line bundles in \cite{DNW19} and on vector bundles in \cite{DNWZ20}.
Theorem \ref{thm:qpositive} is a generalization for partial positivity of the result obtained by the authors in \cite{DNW19}.

%\begin{corollary}\label{cor:rcline}
%	Let $D, \omega, L, h$ and $c$ be the same notation as in Theorem \ref{thm:qpositive}. Then the following properties are equivalent: 
%	\begin{enumerate}
%		\item $\tr_\omega (\ai \Theta_{(L, h)}) \geq c$. 
%		\item For any smooth strictly plurisubharmonic function $\psi$ and any smooth $L$-valued $(n,n)$-form $f$ with compact support, 
		%$$
		%\int_{D} \frac{1}{\tr_\omega (\ddbar\psi) + %c}|f|^2_{(\omega, h)}e^{-\psi}dV_\omega < +\infty,
		%$$
%		there exists $L$-valued $(n,n-1)$-form $u$ satisfying 
%		$\dbar u=f$ and 
%		$$
%		\int_{D}|u|^2_{(\omega, h)}e^{-\psi} \leq \int_{D} \frac{1}{\tr_\omega (\ddbar\psi) + c}|f|^2_{(\omega, h)}e^{-\psi}dV_\omega.
%		$$ 
%	\end{enumerate}
%\end{corollary}

%\begin{theorem}\label{thm:rcvector}
		%Let $D$ be a bounded pseudoconvex domain in $\mathbb{C}^n$, $\omega$ be a K\"ahler form on $D$, and $E\to D$ be a trivial vector bundle over $D$. We also let $h$ be a smooth Hermian metric on $E$, $a$ be a non-zero holomorphic section of $E$, and $c\geq 0$ be some semi-positive constant. Then the following properties are equivalent: 
	%	\begin{enumerate}
		%	\item $\tr_\omega((\ai \Theta_{(E, h)}a, a)_h) \geq c$
			%\item For any smooth strictly plurisubharmonic function $\psi$ and any $E$-valued $(n,n)$-form with finite $L^2$-norm 
		%	$$
		%	\int_{D} 
			%$$
	%	\end{enumerate}
%\end{theorem}

As a higher-rank analogue, we also establish a characterization of RC-positivity. 
RC-positivity is a partial positivity notion introduced by Yang in \cite{Yan18} and a higher-rank analogue of $(\dim X-1)$-positivity.  
We use the same notation as in Theorem \ref{thm:qpositive}.
%Using the following theorem, we show the RC-positivity of a certain kind of direct image bundle (= Theorem \ref{thm:directrc}). 

\begin{theorem}\label{thm:rc}
	Let $E\to D$ be a  vector bundle %over a bounded domain $D\subset \mathbb{C}^n$, $\omega$ be the standard K\"ahler metric on $D$, 
	and $h$ be a smooth Hermitian metric on $E$. % and $c\geq 0$ be some non-negative constant. 
	%Suppose that $E$ is finite rank or a trivial  vector bundle of infinite rank. 
	%In the latter case, the fiber is a separable Hilbert space. 
	Assume that %the following condition: 
	if $\psi$ is a smooth strictly plurisubharmonic function on $D$ and $f$ is a smooth $E$-valued $(n,n)$-form with compact support, there exists a solution of $\dbar u=f$  satisfying
	$$
	\int_D |u|^2_{(\omega, h)} e^{-\psi}dV_\omega\leq \int_D \langle ([\ddbar \psi \otimes Id_E, \Lambda_\omega]+c)^{-1}f, f\rangle_{(\omega, h)}e^{-\psi}dV_\omega.
	$$ 
	Then we obtain 
	$$
	\tr_\omega(\ai\Theta_{(E, h)}a, a)_h(x)\geq c|a|^2_h(x)
	$$
	for any point $x\in D$ and  any element $a\in E_x$. 
	Especially, $(E,h)$ is RC-positive if $c>0$ and 
	RC-semi-positive if $c=0$. 
\end{theorem}

As an application of the characterization in Theorem \ref{thm:qpositive}, 
we prove that uniform $q$-positivity is preserved with respect to a decreasing sequence (Theorem \ref{lem:decreasingq}).
This property is well-known in the case that $q=0$, that is, it is a sequence of plurisubharmonic functions. 
We also propose the definition of uniform $q$-positivity for singular Hermitian metrics (Definition \ref{def:shmpartial}).
%we generalize the  Pr\'ekopa-Berndtsson theorem. 
%In \cite{Pre73}, Pr\'ekopa proved the following theorem. 
%If $\phi(t, x)$ is a convex function on $\mathbb{R}^r_t\times \mathbb{R}^n_x$, 
%then the function $\widetilde{\phi}(t)$ on $\mathbb{R}^r$ defined by 
%$$
%e^{-\widetilde{\phi}(t)}=\int_{x\in \mathbb{R}^n} e^{-\phi(t, x)} dx
%$$
%is convex. 

%Replacing a convex function by a plurisubharmonic function, we can consider a complex version of Pr\'ekopa's theorem. 
%However, it is known that the complex version of Pr\'ekopa's theorem is not true without any additional assumptions \cite{Kis78}. 
%This counterexample is given by Kiselman in \cite{Kis78}. 
%In \cite{Kis78}, Kiselman showed the minimum principle by assuming that the plurisubharmonic function satisfies some invariant properties. 
%Kiselman's minimum principle can be regarded as a special case of the complex Pr\'ekopa theorem, which was mentioned in \cite{Ber98}. 
%In \cite{Ber98}, Berndtsson proved the complex version of Pr\'ekopa's theorem  by assuming that the plurisubharmonic function satisfies some invariant properties. 
%In this article, we call the following result the Pr\'ekopa-Berndtsson theorem.

As a further study, we propose the following problem, which generalizes the Pr\'ekopa-Berndtsson theorem (\cite[Theorem 1.3]{Ber98}).

\begin{problem}\label{thm:qprekopa}
	Let $U$ be a bounded domain in $\mathbb{C}_{z}^n$ and $D$ be a bounded pseudoconvex domain in $\mathbb{C}_{w}^m$. 
	Let $\varphi$ be a smooth function on $\overline{U_z\times D_w}\subset \mathbb{C}_{z}^n \times \mathbb{C}_{w}^m$. 
	We set $\omega_0, \omega_1$ and $\omega_2$ be the standard K\"ahler metrics on $U\subset \mathbb{C}_{z}^n$, $D\subset \mathbb{C}_{w}^m$ and $U\times D\subset \mathbb{C}_{z}^n \times \mathbb{C}_{w}^m$, respectively.
	Assume that 
	\begin{enumerate}
		\item $D$ is a connected Reinhardt domain and $\varphi(z,w_1, \cdots,w_m)$ is independent of $\arg (w_j)$ 
		for $1\leq j\leq m$.
		\item The summation of any distinct $q$ eigenvalues of $\ddbar \varphi$ with respect to $\omega_2$ is greater than or equal to $c$, where $c\geq 0$ is a non-negative constant and $1\leq q\leq n$. 
	\end{enumerate}
We define the function $\phinyoro$ on $U$ by 
$$
e^{-\phinyoro(z)}:=\int_{w\in D}e^{-\varphi(z, w)}d\omega_1(w).
$$
Then the summation of any distinct $q$ eigenvalues of $\ddbar\phinyoro$ 
with respect to $\omega_0$ is greater than or equal to $c$. 
\end{problem}

We immediately see that Problem \ref{thm:qprekopa} is true in the case that $q=1$, which corresponds to the Pr\'ekopa-Berndtsson theorem. 
The proof of the Pr\'ekopa-Berndtsson theorem is based on H\"ormander's $L^2$-estimates, and on a partial converse of these estimates in one variable. 
In Section \ref{sec:shm}, we explain the reason why Theorem \ref{thm:qpositive} can be applied to Problem \ref{thm:qprekopa}. 

The organization of this paper is as follows. 
In Section \ref{sec:prelimi}, we introduce definitions of $q$-positivity, uniform $q$-positivity and RC-positivity. 
We also explain the result of H\"ormander's $L^2$-estimate which we use in this article.
In Section \ref{sec:qpositivity}, we characterize partial positivity by using the H\"ormander $L^2$-estimate. We also show the proofs of Theorem \ref{thm:qpositive} and \ref{thm:rc}.
In Section \ref{sec:appli}, we show applications of the main theorems. We also discuss a definition of uniform $q$-positivity for singular Hermitian metrics. 
%In Section \ref{sec:prekopa}, we prove Theorem \ref{thm:qprekopa}.
%In Section \ref{sec:directrc}, we discuss the RC-positivity of the vector bundle in Theorem \ref{thm:directrc} and propose a conjecture.
In Section \ref{sec:shm}, we propose some problems. 

\vskip10mm
{\bf Acknowledgment. }
The author would like to thank his supervisor Prof. Shigeharu Takayama for enormous supports. 
The author specially wishes to express his gratitude to Dr. Wang Xu for his pointing out mistakes in the first version of the manuscript. 
He would also be grateful to Prof. Shin-ichi Matsumura for helpful comments. 
This work is supported by the Program for Leading Graduate Schools, MEXT, Japan. 
This work is also supported by JSPS KAKENHI Grant Number 18J22119. 

\section{Preliminaries}\label{sec:prelimi}
%Throughout this paper, we use the following notations. 
%\begin{notation}

\textbf{Notation.}
\begin{itemize}
%\item $K_X$ : the canonical line bundle of $X$.
%\item $(X, \omega_X)$ : 
\item $dV_\omega:=\frac{\omega^n}{n!}$ : the volume form determined by $\omega$. 
%\item $E^\star$ : the dual bundle of $E$.
%\item $h^\star$ : the dual metric of $h$ on $E^\star$.
%\item $\mathscr{O}(E)$ : the sheaf of germs of local holomorphic sections of $E$.
\item $C^k_{(p, q)}(X, E):=C^k(X, \wedge^{p, q}T^\star_X\otimes E)$ for $0\leq k\leq +\infty$. 
\item $\mathscr{D}_{(p,q)}(X, E)$ : the space of smooth sections of $\wedge^{(p, q)}T^\star_X\otimes E$ with compact support. 
\item $L^2_{(p, q)}(X, E ; \omega, h)$ : the space of $L^2$ sections of $\wedge^{p, q}T^\star_X\otimes E$ with respect to $\omega$ and $h$. 
\item $\llangle \alpha, \beta \rrangle_{(\omega, h)} := \int_X \langle\alpha, \beta\rangle_{(\omega, h)}dV_\omega$.
\item $\| \alpha\|^2_{(\omega, h)}:= \llangle \alpha, \alpha \rrangle_{(\omega, h)}$.
\item $D'^\star_\psi$ : the adjoint operator of $D'_{\psi}$ with respect to $\llangle \cdot,\cdot  \rrangle_{(\omega, he^{-\psi})}$.
\item $\dbar^\star_\psi$ : the adjoint operator of $\dbar$ with respect to $\llangle \cdot,\cdot  \rrangle_{(\omega, he^{-\psi})}$.
\item $\Delta'_\psi:=D'_\psi D'^\star_\psi+D'^\star_\psi D'_\psi, \Delta_\psi''=\dbar \dbar^\star_\psi + \dbar^\star_\psi \dbar$ with respect to $\llangle \cdot,\cdot  \rrangle_{(\omega, he^{-\psi})}$. 
\item $L_\omega$ : the operator defined by $\omega \wedge \cdot $.
\item $\Lambda_\omega$ : the adjoint operator of $L_\omega$. 
\item $[\cdot, \cdot]$ : the graded Lie bracket. 
%\item $\Delta^n(p ; r):= \{ (z_1, \cdots z_n)\in \mathbb{C}^n \mid |z_i - p_i|< r \}$ for $p= (p_1, \cdots , p_n)\in \mathbb{C}^n$. 
%\item $\Delta^n_r:= \Delta^n(o; r)$. 
\item $\mathbb{B}^n_r:= \{ (z_1, \cdots, z_n)\in\mathbb{C}^n \mid \sum_{i=1}^{n}|z_i|^2<r^2 \}$.
\end{itemize}
%\end{notation}

In \cite{AG62}, Andreotti and Grauert introduced partial positivity notions and studied partially vanishing cohomology. 
Here, we introduce the notions of $q$-positivity and uniform $q$-positivity for smooth Hermitian metrics on line bundles. 
\begin{definition}(cf. \cite{AG62}, \cite[Definition 2.1]{Yan19})\label{def:qpositive}
	Let $L\to X$ be a holomorphic line bundle over a complex manifold $X$ with $\dim X=n$. 
	Let $h$ be a smooth Hermitian metric on $L$. 
	We say that 
	\begin{enumerate}
		\item $ (L,h) $ is {\it $ q $-$($semi-$)$positive} if the Chern curvature $\ai\Theta_{(L, h)}$ has at least $(n-q)$ (semi-)positive eigenvalues at any point on $X$. We also say that $L$ is {\it $ q $-$($semi-$)$positive} if there exists a smooth Hermitian metric $h$ on $L$ such that $(L,h)$ is $q$-(semi-)positive. 
		\item $ (L,h) $ is {\it uniformly $ q $-$($semi-$)$positive} if there exists a smooth Hermitian metric $\omega$ such that the summation of any distinct $(q+1)$ eigenvalues of the Chern curvature $\ai\Theta_{(L, h)}$ with respect to $\omega$ is (semi-)positive at any point on $X$. 
		We also say that $L$ is {\it uniformly $q$-$($semi-$)$positive} if there exist a smooth Hermitian metric $h$ on $L$ and a smooth Hermitian metric $\omega$ such that $(L, h)$ is uniformly $q$-(semi-)positive with respect to $\omega$. 
	\end{enumerate}
\end{definition}

A simple computation yields that uniform $q$-(semi-)positivity implies $q$-(semi-)positivity. 
Note that usual (semi-)positivity corresponds to $0$-(semi-)positivity. 
Conversely, it is known that the above two positivity notions are equivalent over a compact complex manifold. 

\begin{proposition}$($\cite[Proposition 2.2]{Yan19}$)$\label{prop:uniform}
	Let $X$ be a compact complex manifold and $L$ be a $q$-positive line bundle. Then $L$ is a uniformly $q$-positive line bundle.
\end{proposition}

Next, we also give definitions of RC-positivity and weak RC-positivity, which were introduced by Yang in \cite{Yan18}.

\begin{definition}(\cite[Definition 3.3]{Yan18})\label{def:rcpositivity}
	A Hermitian holomorphic vector bundle $(E, h)$ over a complex manifold $X$ is called {\it RC-positive} (resp. {\it RC-negative}) if at any point $x\in X$ and for any non-zero element $a\in E_x$, there exists a vector $v\in T_xX$ such  that 
	$$
	(\ai \Theta_{(E, h)}(v,v)a, a)_h >0 ~( \text{resp.} <0).
	$$
	We also call $(E,h)$ {\it weakly RC-positive} if there exists a smooth Hermitian metric $h$ on the tautological line bundle $\mathscr{O}_E(1)$ over $\mathbb{P}(E^\star)$ such that $(\mathscr{O}_E(1), h)$ is ($\dim X-1$)-positive. 
\end{definition}

Note that Griffiths positivity implies RC-positivity. 
Moreover, if $\dim X=1$, RC-positivity is equivalent to Griffiths positivity. If ${\rm rank} E=1$, RC-positivity is the same concept as $(\dim X-1)$-positivity.

Finally, we mention the following result, which was initially obtained by H\"ormander \cite{Hor65}.
H\"ormander's $L^2$-estimate is fundamental  and important in several complex variables. 
In our paper, we use this $L^2$-estimate to characterize several notions of partial positivity. 
Here, we adopt the following form.

\begin{theorem}$($cf. \cite{Dem82}, \cite[Theorem (5.1)]{Dem12} and \cite[Theorem 6.1 in Chapter VIII]{DemCom}$)$\label{thm:demailly}
	Let $(X, \widehat{\omega})$ be a complete K\"ahler manifold, $\omega$ be another K\"ahler metric which is not necessarily complete, and $(E, h)\to X$ be a holomorphic line bundle. 
	We also let $A_{(\omega, h)}=[\sqrt{-1}\Theta_{(E, h)}, \Lambda_\omega]$ be the curvature operator in bidegree $(n, q)$ for $q\geq 1$. 
	Assume that $A_{(\omega, h)}$ is positive definite everywhere on $\wedge^{n,q}T^\star X\otimes E$. 
	Then for any $\dbar$-closed $f\in L^2_{(n, q)}(X, E ; \omega, h)$, there exists $u\in L^2_{(n, q-1)}(X, E ; \omega, h)$
	such that $\dbar u=f$ and 
	$$
	\int_X |u|^2_{(\omega, h)}dV_\omega \leq \int_X \langle A_{(\omega, h)}^{-1}f, f\rangle_{(\omega, h)}dV_\omega, 
	$$
	where we assume that the right-hand side is finite. 
\end{theorem}

%\begin{proof}
%	As this is a known result, we give a sketch of the proof. 
%	If $\omega$ is complete, the theorem is a usual H\"ormander's $L^2$-estimates. If $\omega$ is not complete, we set the K\"ahler metric 
%	$$
%	\omega_\epsilon = \omega + \epsilon \widehat{\omega}
%	$$
%	for $\epsilon >0$, which is complete.
%	It follows that 
%	$$
%	|f|^2_{(\omega_\epsilon,h)}dV_{\omega_\epsilon}\leq |f|^2_{(\omega,h)}dV_\omega, \hspace{5mm} \langle A_{(\omega_\epsilon, h)}^{-1}f, f\rangle_{(\omega_\epsilon, h)}dV_{\omega_\epsilon}\leq \langle A_{(\omega, h)}^{-1}f, f \rangle_{(\omega, h)}dV_\omega
%	$$ 
%	for any $\dbar$-closed $f\in L^2_{(n,q)}(X, E;\omega,h)$. Then we get a solution $u_\epsilon\in L^2(X, E;\omega,h)$ such that 
%	$\dbar u_\epsilon=f$ and 
%	$$
%	\int_X |f_\epsilon|^2_{(\omega_\epsilon,h)}dV_{\omega_\epsilon}
%	\leq \int_X \langle A_{(\omega_\epsilon, h)}^{-1}f, f\rangle_{(\omega_\epsilon, h)}dV_{\omega_\epsilon}
%	\leq \int_X \langle A_{(\omega, h)}^{-1}f, f\rangle_{(\omega, h)}dV_{\omega}. 
%	$$
%	Finally, it is enough to take the weak limit $f$ of $(f_\epsilon)$. 
%\end{proof}

\section{A characterization of partial positivity via $L^2$-estimates}\label{sec:qpositivity}
\subsection{Uniform $q$-positivity}

In this subsection, we discuss a characterization of uniform $q$-positivity in terms of $L^2$-estimates.
In order to prove Theorem \ref{thm:qpositive}, 
we need the following lemma. 
The proof is a simple computation. 

\begin{lemma}$($cf. \cite[(4.10)]{Dem12} and \cite[Proposition (5.8) in Chapter VI]{DemCom}$)$\label{lem:nq}
	Let the notation be the same as in Theorem \ref{thm:qpositive}. We also let $f$ be any $\dbar$-closed $L$-valued $(n,q)$-form.
	%$$
	%f=\sum_{1\leq i_1<\cdots<i_q\leq n}f_{i_1\cdots i_q}dz_1\wedge \cdots \wedge dz_n\wedge d\bar{z}_{i_1}\wedge \cdots \wedge d\bar{z}_{i_q} \otimes e_L.
	%$$
	At a fixed point $p\in X$, we take a coordinate $(z_1, \cdots , z_n)$ around $p$ such that 
	$$
	\omega=\ai \sum_{j=1}^{n} dz_j\wedge d\bar{z}_j, \hspace{10mm}\ai\Theta_{(L, h)}= \ai \sum_{j=1}^{n}\gamma_j dz_j\wedge d\bar{z}_j.
	$$
	We write 
	$$
	f=\sum_{1\leq i_1<\cdots<i_q\leq n}f_{i_1\cdots i_q}dz_1\wedge \cdots \wedge dz_n\wedge d\bar{z}_{i_1}\wedge \cdots \wedge d\bar{z}_{i_q} \otimes e_L
	$$
	for a local holomorphic frame $e_L$ of $L$ around $p$. 
	Then we get 
	$$
	[\ai\Theta_{(L, h)}, \Lambda_\omega]f = \sum_{1\leq i_1<\cdots<i_q\leq n}(\sum_{k=1}^{q} {\gamma_{i_k}})f_{i_1\cdots i_q}dz_1\wedge \cdots \wedge dz_n\wedge d\bar{z}_{i_1}\wedge \cdots \wedge d\bar{z}_{i_q} \otimes e_L.
	$$
\end{lemma}

We now give a proof of Theorem \ref{thm:qpositive}.
The idea for the proof is based on the arguments in \cite[Theorem 2.1]{DNW19} and \cite[Theorem 3.1]{DNWZ20}.

\begin{proof}[\indent\sc Proof of Theorem \ref{thm:qpositive}]\label{prf:rcline}
	$(1)\Longrightarrow (2)$. 
	We have 
	$$
	[\ai \Theta_{(L, he^{-\psi})}, \Lambda_\omega] = [\ai \Theta_{(L, h)}, \Lambda_\omega] + [\ddbar \psi, \Lambda_\omega]
	$$
	for any smooth strictly plurisubharmonic function $\psi$. 
	We fix a smooth $\dbar$-closed $L$-valued $(n, q)$-form $f$ with compact support. 
	%$$
	%\int_{D} \langle ([\ddbar \psi, \Lambda_\omega]+c)^{-1}f, %f\rangle_{(\omega, h)} e^{-\psi}dV_\omega < +\infty.
	%$$ 
	%We can compute the curvature operator $A_{(\omega, h)}$ as follows
	%\begin{align*}
	%\langle [\ai \Theta_{(L, h)}, \Lambda_\omega]f, f\rangle_{(\omega, h)} &=\tr_\omega(\ai \Theta_{(L, h)})|f|^2_{(\omega, h)}.
	%\langle [\ddbar \psi, \Lambda_\omega]f, f\rangle_{(\omega, h)} &=\tr_\omega(\ddbar \psi)|f|^2_{(\omega, h)}.
	%\end{align*}
	The assumption of $(1)$ and Lemma \ref{lem:nq} implies that
	% $\tr_\omega(\ai \Theta_{(L, h)}) \geq c \geq 0$. 
	\begin{align*}
	\langle [\ai \Theta_{(L, h)}, \Lambda_\omega]f, f\rangle_{(\omega, h)} &\geq c|f|^2_{(\omega, h)}.
	\end{align*}
	The curvature operator $[\ai \Theta_{(L, he^{-\psi})}, \Lambda_\omega]$ is positive definite on $\wedge^{n,q}T^\star D\otimes L$ everywhere. 
	 Therefore, by using Theorem \ref{thm:demailly}, we can solve the $\dbar$-equation $\dbar u =f$ as follows 

\begin{align*}
\int_{D}|u|^2_{(\omega,h)}e^{-\psi}dV_\omega &\leq 	\int_{D} \langle [\ai \Theta_{(L, he^{-\psi})}, \Lambda_\omega]^{-1} f, f\rangle_{(\omega, h)}e^{-\psi}dV_\omega\\
& \leq \int_{D} \langle ([\ddbar \psi, \Lambda_\omega]+c)^{-1}f, f\rangle_{(\omega, h)} e^{-\psi}dV_\omega < +\infty
\end{align*}
	for some $u\in L^2_{(n, q-1)}(D, L ; \omega, he^{-\psi})$.

$(2)\Longrightarrow (1)$. For any smooth strictly plurisubharmonic function $\psi$ and any $\dbar$-closed $f\in \mathscr{D}_{(n, q)}(D, L)$, %with finite $L^2$-norm 
%$$
%\int_{D} \langle ([\ddbar \psi, \Lambda_\omega]+c)^{-1}f, f\rangle_{(\omega, h)}e^{-\psi}dV_\omega < +\infty, 
%$$
we get a solution $u\in L^2_{(n, q-1)}(D, L ; \omega, he^{-\psi})$ of $\dbar u =f$ satisfying 
$$
\| u\|^2_{(\omega, he^{-\psi})}\leq \llangle ([\ddbar \psi, \Lambda_\omega]+c)^{-1}f,  f\rrangle_{(\omega, he^{-\psi})}.
$$
Set $g:= ([\ddbar \psi, \Lambda_\omega]+c)^{-1}f$. We obtain 
\begin{align*}
|\llangle g, f\rrangle_{(\omega, he^{-\psi})}|^2&=|\llangle g, \dbar u \rrangle_{(\omega, he^{-\psi})}|^2\\
&\leq |\llangle \dbar^\star_\psi g, u\rrangle_{(\omega, he^{-\psi})}|^2\\
&\leq \| \dbar^\star_\psi g\|^2_{(\omega, he^{-\psi})} \| u\|^2_{(\omega, he^{-\psi})}\\
&\leq \| \dbar^\star_\psi g\|^2_{(\omega, he^{-\psi})} |\llangle g, f \rrangle_{(\omega, he^{-\psi})}|. 
\end{align*}
Using the Bochner-Kodaira-Nakano identity $\Delta''_\psi=\Delta_\psi'+[\ai \Theta_{(L, he^{-\psi})}, \Lambda_\omega]$ (cf. \cite[(4.6)]{Dem12}), we have 
\begin{align*}
\| \dbar^\star_\psi g\|^2_{(\omega, he^{-\psi})} &= \llangle (\Delta_\psi''-\dbar_\psi^\star \dbar)g,  g\rrangle_{(\omega, he^{-\psi})}\\
&= \llangle \Delta_\psi'g, g \rrangle_{(\omega, he^{-\psi})} + \llangle [\ai\Theta_{(L, he^{-\psi})}, \Lambda_\omega]g, g\rrangle_{(\omega, he^{-\psi})} -\| \dbar g \|^2_{(\omega,he^{-\psi})} \\
&\leq  \| D'^\star_\psi g \|^2_{(\omega, he^{-\psi})} + \llangle [\ai\Theta_{(L, h)}, \Lambda_\omega] g, g\rrangle_{(\omega, he^{-\psi})} + \llangle [\ddbar \psi, \Lambda_\omega]g, g \rrangle_{(\omega, he^{-\psi})}.
\end{align*}
Therefore, we get 
\begin{align*}
&\llangle g, ([\ddbar \psi, \Lambda_\omega]+c)g\rrangle_{(\omega, he^{-\psi})}\\
&\leq \| D'^\star_\psi g \|^2_{(\omega, he^{-\psi})} + \llangle [\ai\Theta_{(L, h)}, \Lambda_\omega] g, g\rrangle_{(\omega, he^{-\psi})} + \llangle [\ddbar \psi, \Lambda_\omega]g, g \rrangle_{(\omega, he^{-\psi})},
\end{align*}
that is, 
\begin{equation}\label{eq:sekibun}
0 \leq \| D'^\star_\psi g \|^2_{(\omega, he^{-\psi})} + \llangle ([\ai\Theta_{(L, h)},\Lambda_\omega]-c)g, g\rrangle_{(\omega, he^{-\psi})}.
\end{equation}

We give a proof by contradiction.
In other words, we suppose that the summation of some distinct $q$ eigenvalues of $\ai\Theta_{(L, h)}$ with respect to $\omega$ is less than $c$ at some point $a\in D$. 
We can assume that $o=a\in D$, where $o$ is the origin of $\mathbb{C}^n$.
Let $\gamma_1, \cdots, \gamma_n$ be the eigenvalues of $\ai\Theta_{(L, h)}$ with respect to $\omega$, which are globally defined on $D$. 
Changing the coordinate by some unitary transformation, we take a coordinate $(z_1, \cdots, z_n)$ centered at $o$ such that 
$$
\omega=\ai \sum dz_j\wedge d\bar{z}_j,
$$
on $D$ and 
$$
\ai \Theta_{(L, h)}=\ai \sum \gamma_j dz_j\wedge d\bar{z}_j
$$
at $o$. 
Without any loss of generality, we suppose that 
$$
\gamma_1(o) +\cdots +\gamma_q(o) - c<0.
$$
We fix an open neighborhood $U$ of $o$ and a local holomorphic frame $e_L$ of $L$ on $U$. 
We define 
$$
F:=dz_1\wedge \cdots\wedge dz_n\wedge d\bar{z}_1\wedge \cdots\wedge d\bar{z}_q\otimes e_L \hspace{3mm}\in C^\infty_{(n,q)}(U, L).
$$
Then we have 
\begin{align*}
&\langle ([\ai\Theta_{(L, h)}, \Lambda_\omega]-c)F, F\rangle_{(\omega, h)}(o)\\
&=\langle (\sum_{j=1}^{q}\gamma_j-c)F, F\rangle_{(\omega, h)}(o)\\
&=(\sum_{j=1}^{q}\gamma_j(o)-c)|e_L|^2_h<0.
\end{align*}
We take a positive constant $\delta>0$ such that 
$$
\langle ([\ai\Theta_{(L, h)}, \Lambda_\omega]-c)F, F\rangle_{(\omega, h)}(o)=(\sum_{j=1}^{q}\gamma_j(o)-c)|e_L|^2_h = -2\delta. 
$$
Since $\langle ([\ai\Theta_{(L, h)}, \Lambda_\omega]-c)F, F\rangle_{(\omega, h)}$ has continuous coefficients, 
we take a sufficiently small $r\in (0, +\infty)$ such that 
$\mathbb{B}^n_r \subset U\Subset D$ and 
$$
\langle ([\ai\Theta_{(L, h)}, \Lambda_\omega]-c)F, F\rangle_{(\omega, h)} < -\delta
$$
on $\mathbb{B}^n_r$. 

We take a smooth strictly plurisubharmonic function $\psi(z)=|z|^2-\frac{r^2}{4}$ on $D$. 
Let $\chi$ be a cut-off function on $\mathbb{B}^n_r$ such that 
$\chi$ is smooth, $0\leq \chi \leq 1$, supp$\chi\Subset \mathbb{B}^n_r$ and $\chi|_{\mathbb{B}^n_{r/2}}\equiv 1$. 
We set $v:= (-1)^{n+q-1}\chi \bar{z}_q dz_1\wedge \cdots \wedge dz_n\wedge d\bar{z}_1\wedge \cdots \wedge d\bar{z}_{q-1}\otimes e_L$ and $g:=\dbar v$. 
Then $g$ is a $\dbar$-closed $L$-valued $(n,q)$-form with compact support and 
$$
g=dz_1\wedge \cdots \wedge dz_n\wedge d\bar{z}_1\wedge \cdots \wedge d\bar{z}_{q}\otimes e_L
$$
on $\mathbb{B}^n_{r/2}$.
We remark that $[\ddbar(m\psi), \Lambda_\omega]f=mqf$ for $f\in \wedge^{n,q}T^\star_D\otimes L$.
We define $f_m:= ([\ddbar(m\psi), \Lambda_\omega]+c)g=(mq+c)g$. 
%$f_m:=(\tr_\omega(\ddbar(m\psi))+c)g=(mn+c)g$ for $m\in \mathbb{N}$. 
%We regard $g$ and  $f_m$ as $L$-valued $(n,n)$-forms on $D$.
It clearly holds that $f_m$ is an also $\dbar$-closed $L$-valued $(n,q)$-form with compact support. 
%Since $g$ has compact support, it follows that $f_m\in L^2(D, L;\omega, he^{-m\psi})$ and 
%$$
%\int_{D} \langle ([\ddbar(m\psi), \Lambda_\omega]+c)^{-1}f_m, f_m\rangle_{(\omega, h)} e^{-m\psi}dV_\omega
%=\int_{D} \langle g, f_m\rangle_{(\omega, h)} e^{-m\psi}dV_\omega
%< +\infty
%$$
%for every $m\in \mathbb{N}$. 
%Set $g_m:= (
%[\ddbar(m\psi), \Lambda_\omega]+c)^{-1}f$. 
Then $g$ satisfies the inequality (\ref{eq:sekibun}) for every $m\psi$.
Considering the commutation relation $\ai[\Lambda_\omega,\dbar]=D'^\star_{m\psi}$ (cf. \cite[(1.1) in Chapter VII]{DemCom}), we have that 
$$
D'^\star_{m\psi}g= 0
$$
on $\mathbb{B}^n_{r/2}$ since $\omega$ is the standard K\"ahler metric and $g$ has constant coefficients on $\mathbb{B}^n_{r/2}$, and 
$$
|D'^\star_{m\psi}g|^2_{(\omega, h)}\leq C_1
$$ 
for some positive constant $C_1>0$ which is independent of $m$ and $\psi$ on $\mathbb{B}^n_r$. 

Since $g=F$ on $\mathbb{B}^n_{r/2}$, 
we know that $\langle ([\ai\Theta_{(L, h)}, \Lambda_\omega]-c)g, g\rangle_{(\omega, h)} < -\delta$ on $\mathbb{B}^n_{r/2}$ and $\langle ([\ai\Theta_{(L, h)}, \Lambda_\omega]-c)g, g\rangle_{(\omega, h)}\leq C_2 $ for some positive constant $C_2>0$ on $\mathbb{B}^n_r$. 
%Let $C_3:=\sup_{\mathbb{B}^n_r}|g|^2_{(\omega,h)}<+\infty$. 
Consequently, we can compute the right-hand side of (\ref{eq:sekibun}) for $g$ and $m\psi$ as follows: 
\begin{align*}
0\leq& \int_{D}|D'^\star_{m\psi}g|^2_{(\omega,h)}e^{-m\psi}dV_\omega + \int_{D} \langle ([\ai\Theta_{(L, h)}, \Lambda_\omega]-c)g, g\rangle_{(\omega, h)}e^{-m\psi}dV_\omega\\
=&\int_{\mathbb{B}^n_r\setminus \overline{\mathbb{B}^n_{r/2}}}|D'^\star_{m\psi}g|^2_{(\omega,h)}e^{-m\psi}dV_\omega + \int_{\mathbb{B}^n_{r/2}}\langle ([\ai\Theta_{(L, h)}, \Lambda_\omega]-c)g, g\rangle_{(\omega, h)}e^{-m\psi}dV_\omega\\
&+\int_{\mathbb{B}^n_r\setminus \overline{\mathbb{B}^n_{r/2}}}\langle ([\ai\Theta_{(L, h)}, \Lambda_\omega]-c)g, g\rangle_{(\omega, h)}e^{-m\psi}dV_\omega\\
\leq& (C_1+C_2)\int_{\mathbb{B}^n_r\setminus \overline{\mathbb{B}^n_{r/2}}}e^{-m\psi}dV_\omega -\delta \int_{\mathbb{B}^n_{r/2}}e^{-m\psi}dV_\omega.
\end{align*}
Since $\psi > 0$ on $ \mathbb{B}^n_r\setminus \overline{\mathbb{B}^n_{r/2}} $, the first term goes to zero as $m\to +\infty$ by Lebesgue's dominated convergence theorem. 
The second term has a negative upper bound 
$$
-\delta  \int_{\mathbb{B}^n_{r/2}}e^{-m\psi}dV_\omega < -\delta  |\mathbb{B}^n_{r/2}|
$$
which is independent of $m$ since $\psi < 0$ on $\mathbb{B}^n_{r/2}$.
Taking a sufficiently large $m>>1$, we get 
$$
(C_1+C_2)\int_{\mathbb{B}^n_r\setminus \overline{\mathbb{B}^n_{r/2}}}e^{-m\psi}dV_\omega -\delta  \int_{\mathbb{B}^n_{r/2}}e^{-m\psi}dV_\omega<0, 
$$
which is a contradiction. 
%Therefore, we can conclude that 
%$ \tr_\omega(\ai\Theta_{(L, h)})\geq c $ on $D$. 
\end{proof}

%\begin{remark}
%	We can prove $(1)\Longrightarrow (2)$ without assuming that $D$ is pseudoconvex. Therefore, we can use the characterization (2) in Theorem \ref{thm:qpositive} to prove Theorem \ref{thm:qprekopa}. 
%\end{remark}

\subsection{RC-positivity}
In this subsection, we give a characterization of RC-positivity via $L^2$-estimates. This is a higher-rank analogue of Theorem \ref{thm:qpositive}. Although the proof is almost identical to the proof of Theorem \ref{thm:qpositive}, we present it for the sake of completeness. 
\begin{proof}[\indent\sc Proof of Theorem \ref{thm:rc}]
	We take an arbitrary smooth strictly plurisubharmonic function $\psi$ and an arbitrary $f\in \mathscr{D}_{(n, n)}(D, E)$.
	Repeating the argument in the proof of Theorem \ref{thm:qpositive} (cf. \cite[Theorem 3.1]{DNWZ20} or \cite[Proposition 2.7]{Ina20}), we obtain the following inequality
	\begin{equation}\label{eq:rcpositive}
	0\leq \| D'^\star_\psi g\|^2_{(\omega,he^{-\psi})}+\llangle ([\ai \Theta_{(E, h)}, \Lambda_\omega]-c)g, g \rrangle_{(\omega, he^{-\psi})},
	\end{equation}
	where $g=([\ddbar\psi\otimes Id_E,\Lambda_\omega]+c)^{-1}f$. 
	
	We give a proof by contradiction. We assume that there exists some point $x\in D$ and some element $a\in E_x\setminus \{0
	\}$ such that 
	\begin{equation}\label{eq:mujunrc}
	\tr_\omega (\ai\Theta_{(E, h)}a, a)_h(x)<c|a|^2_h(x).
	\end{equation}
	We may assume that $x=o\in D$. 
	Since $h$ has smooth
	coefficients, we can take a sufficiently small $r\in (0, +\infty)$ such that $\mathbb{B}^n_r\Subset D$, $E|_{\mathbb{B}^n_r}$ is trivial, and 
	\begin{equation}\label{eq:minus}
	\tr_\omega(\ai\Theta_{(E, h)}a, a)_h-c|a|^2_h<-\delta
	\end{equation}
	on $\mathbb{B}^n_r$ for some positive constant $\delta>0$. 
	Here we regard $a$ as a section of $E$ with constant coefficients.
	
	As 
	in the proof of Theorem \ref{thm:qpositive}, 
	we take a smooth strictly plurisubharmonic function $\psi(z)=|z|^2-\frac{r^2}{4}$ and a cut-off function $\chi$ such that supp$\chi \Subset \mathbb{B}^n_r$ and $\chi|_{\mathbb{B}^n_{r/2
	}}\equiv 1$. 
	We consider the following $E$-valued $(n,n)$-form with compact support 
	$$
	g=\chi a dZ\wedge d\bar{Z}
	$$
	on $D$. 
	Here we use the notation 
	$$
	dZ=dz_1\wedge \cdots \wedge dz_n, \hspace{5mm} d\bar{Z}=d\bar{z}_1\wedge \cdots \wedge d\bar{z}_n
	$$
	for simplicity. 
	We also define 
	$$
	f_m:= ([\ddbar (m\psi)\otimes Id_E, \Lambda_\omega]+c)g=(mn+c)g.
	$$
	Note that $f_m\in \mathscr{D}_{(n, n)}(D, E)$.
	Hence, we see that $g$ satisfies the inequality (\ref{eq:rcpositive}) for each $m\psi$. 
	
	We compute the terms $\langle [\ai \Theta_{(E, h)}, \Lambda_\omega](sdZ\wedge d\bar{Z}), sdZ\wedge d\bar{Z}\rangle_{(\omega, h)}$ and $\tr_\omega (\ai\Theta_{(E, h)}s, s)_h$ for any section $s$ of $E$. 
	Note that $sdZ\wedge d\bar{Z}\in C^\infty_{(n,n)}(D, E)$. We write the curvature tensor $\ai \Theta_{(E, h)}$ as 
	$$
	\ai\Theta_{(E, h)}=\sum_{1\leq j, k\leq n}\Theta_{j\bar{k}}dz_j \wedge d\bar{z}_k, 
	$$
	where $\Theta_{j\bar{k}}$ are operators on each $E_t$. 
	Then we get 
	\begin{align*}
	\langle [\ai \Theta_{(E, h)}, \Lambda_\omega](sdZ\wedge d\bar{Z}), sdZ\wedge d\bar{Z}\rangle_{(\omega, h)}&=\langle (\sum_{j=1}^n \Theta_{j\bar{j}}s)dZ\wedge d\bar{Z}, sdZ\wedge d\bar{Z}\rangle_{(\omega, h)}\\
	&=\sum_{j=1}^n (\Theta_{j\bar{j}}s, s)_h
	\end{align*}
	and 
	\begin{align*}
	\tr_\omega (\ai\Theta_{(E, h)}s, s)_h&=\tr_\omega  \left( \sum_{1\leq j, k\leq n}(\Theta_{j\bar{k}}s, s)_hdz_j\wedge d\bar{z}_k\right)\\
	&= \sum_{j=1}^n (\Theta_{j\bar{j}}s, s)_h.
	\end{align*}
	Hence, on $\mathbb{B}^n_{r/2}$, the inequality (\ref{eq:minus}) implies that 
	$$
	\langle ([\ai \Theta_{(E, h)}, \Lambda_\omega]-c)g, g\rangle_{(\omega,h)}<-\delta. 
	$$
	Then, 
	taking a sufficiently large $m>>1$ and 
	repeating the argument in the proof of Theorem \ref{thm:qpositive} again, we conclude that the inequality (\ref{eq:mujunrc}) contradicts the inequality (\ref{eq:rcpositive}), which completes the proof. 
\end{proof}

\section{Applications of a new characterization}\label{sec:appli}

In this section, we give some applications of Theorem \ref{thm:qpositive}. First, we prove the following theorem. 
Here we use the same notation as in Theorem \ref{thm:qpositive}. 
\begin{theorem}\label{lem:decreasingq}
	Let $\varphi$ be a smooth function on $D$. 
	Suppose that there exists a sequence of smooth functions $\{ \varphi_j\}_{j=1}^\infty$ decreasing to $\varphi$ pointwise such that the summation of any distinct $q$ eigenvalues of $\ddbar\varphi_j$ with respect to $\omega$ is greater than or equal to some non-negative constant $c\geq 0$. 
	Then the summation of any distinct $q$ eigenvalues of $\ddbar\varphi$ with respect to $\omega$ is greater than or equal to $c$.
\end{theorem}

It is well-known that Theorem \ref{lem:decreasingq} holds in the case that $q=1$, that is, $\varphi_j$ are plurisubharmonic functions. 

\begin{proof}
	We use the characterization in Theorem \ref{thm:qpositive}.
	Since the result is a local property, we may assume that $D$ is pseudoconvex. 
	It is enough to show that for any smooth strictly plurisubharmonic function $\psi$ and any smooth $\dbar$-closed $(n,q)$-form $f$ with compact support, there exists a solution of $\dbar u=f$ satisfying 
	$$
	\int_D |u|^2_{\omega_0}e^{-(\varphi+\psi)}dV_{\omega_0}\leq \int_D \langle([\ddbar\psi, \Lambda_{\omega_0}]+c)^{-1}f, f \rangle_{\omega_0}e^{-(\varphi+\psi)}dV_{\omega_0}. 
	$$
	The assumption of $\varphi_j$ implies that we get a solution of $\dbar u_j=f$ satisfying
	\begin{align*}
	\int_D |u_j|^2_{\omega_0}e^{-(\varphi_j+\psi)}dV_{\omega_0}&\leq \int_D \langle([\ddbar\psi, \Lambda_{\omega_0}]+c)^{-1}f, f \rangle_{\omega_0}e^{-(\varphi_j+\psi)}dV_{\omega_0}\\
	&\leq \int_D \langle([\ddbar\psi, \Lambda_{\omega_0}]+c)^{-1}f, f \rangle_{\omega_0}e^{-(\varphi+\psi)}dV_{\omega_0}\\
	&<+\infty
	\end{align*}
	for each $j\in \mathbb{N}$. 
	Note that the right-hand side of the above inequality has an 
	upper bound independent of $j$ and $\{ u_k \}_{k\geq j}$ forms a bounded sequence in $L^2_{(n, q-1)}(D, \mathbb{C};\omega_0, e^{-(\varphi_j+\psi)})$.
	Therefore, 
	 we find a weakly convergent subsequence $\{ u_{j_k}\}_{k=1}^\infty$ by using a diagonal argument and monotonicity of $\{\varphi_j\}_{j=1}^\infty$, 
	 which is the standard argument of $L^2$-solutions of $\dbar$. 
	 We have that $\{ u_{j_k}\}_{k=1}^\infty$ weakly converges in $L^2_{(n, q-1)}(D, \mathbb{C};\omega_0, e^{-(\varphi_j+\psi)})$
	 for every $j$ and the weak limit denoted by $u_\infty$ satisfies $\dbar u_\infty=f$ and 
	 $$
	 \int_{D}|u_\infty|^2_{\omega_0}e^{-(\varphi+\psi)}dV_{\omega_0}\leq \int_D
	  \langle([\ddbar\psi, \Lambda_{\omega_0}]+c)^{-1}f, f \rangle_{\omega_0}e^{-(\varphi+\psi)}dV_{\omega_0}
	 $$
	 due to the monotone convergence theorem. 
	 Then we complete the proof. 
\end{proof}

Next, by using the new characterization, we propose the definition of uniform $q$-positivity for singular Hermitian metrics. 
Note that we can consider the condition $(2)$ in Theorem \ref{thm:qpositive} without assuming that $h$ is smooth. 
\begin{definition}\label{def:shmpartial}
	Let $L$ be a holomorphic line bundle over an $n$-dimensional K\"ahler manifold $(X, \omega)$ and $h$ be a singular Hermitian metric on $L$ such that $-\log h$ is upper semi-continuous. 
	Set $1\leq q \leq n$ and $c \geq 0$. 
	We say that {\it $(L, h)$  is uniformly $(q-1)$-$c$-positive with respect to $\omega$} if for any point $x\in X$, there exists an open neighborhood $U$ of $x$ such that for any relatively compact pseudoconvex domain $D$ in $U$, $(L, h)$, $\omega$ and $c$ satisfy the condition $(2)$ in Theorem \ref{thm:qpositive} on $D$. 
\end{definition}

Thanks to Theorem \ref{thm:qpositive}, in the case that $h$ is smooth, the above definition is equivalent to uniform $(q-1)$-positivity. 
Under this formulation, we can show Theorem \ref{lem:decreasingq} without assuming the condition that $\varphi$ is smooth. 
The proof remains the same. 

\begin{theorem}$($cf. Theorem \ref{lem:decreasingq}$)$
	Let $\{\varphi_j \}_{j=1}^\infty$ be a sequence of smooth functions decreasing to a locally integrable function $\varphi \not\equiv - \infty$.
	If the summation of any distinct $q$ eigenvalues of $\ddbar \varphi_j$ with respect to $\omega$ is greater than or equal to $c$,  $(\mathbb{C}, e^{-\varphi})$ is uniformly $(q-1)$-$c$-positive in the sense of Definition \ref{def:shmpartial}. 
\end{theorem}

The argument above has many other applications. 
For instance, it is known that Nakano semi-positivity can be characterized via $L^2$-estimates (cf. \cite[Theorem 1.1]{DNWZ20}). 
By using the same method, we can also show that if a sequence of smooth Nakano semi-positive Hermitian metrics $\{h_\nu\}_{\nu=1}^\infty$ increasing to a (possibly singular) Hermitian metric $h$, $h$ is also Nakano semi-positive (for Nakano semi-positivity of singular Hermitian metrics, see \cite[Definition 1.2]{Ina20}).

\section{Further study}\label{sec:shm}

In this section, we propose some problems. 
First, we discuss Problem \ref{thm:qprekopa}.
In \cite{Ber98}, Berndtsson generalized the Pr\'ekopa theorem \cite{Pre73} by assuming that the plurisubharmonic function satisfies some invariant properties. 
\begin{theorem}$($\cite[Theorem 1.3]{Ber98}$)$\label{thm:prekopa}
	Let $\varphi $ be a plurisubharmonic function on $U_{z}\times D_{w}\subset \mathbb{C}_{z}^n\times \mathbb{C}_{w}^m$, where $D_w$ is pseudoconvex. 
	Assume that one of the following conditions holds:
	\begin{enumerate}
		\item $D$ is a connected Reinhardt domain and $\varphi(z, w_1, \cdots , w_m)$ is independent of $\arg (z_j)$ for $1\leq j \leq m$. 
	\item $D$ contains the origin and for any $z\in U$, $w\in D$, and $\theta \in \mathbb{R}$, we have 
	$e^{\ai \theta}w\in D$ and $\varphi(z, e^{\ai \theta}w)=\varphi(z, w)$. 
	\end{enumerate}
	Then the function $\phinyoro$ defined on $U$ by 
$$
e^{-\phinyoro(z)}:= \int_{w\in D}e^{-\varphi(z, w)}
$$
is plurisubharmonic. 
\end{theorem} 

This research has been generalized in a variety of directions (cf. \cite{Cor05}, \cite{DZZ14}).
Since Theorem \ref{thm:prekopa} can be applied in the case where $D$ is bounded, Problem \ref{thm:qprekopa} is a generalization of Theorem \ref{thm:prekopa} for partial positivity. 
Here we explain the reason why Theorem \ref{thm:qpositive} is one strategy to prove Problem \ref{thm:qprekopa}. 
Let $\pi:U\times D\to U$ be the projection map and $dZ=dz_1\wedge \cdots \wedge dz_n, dW=dw_1\wedge \cdots \wedge dw_m$. 

Thanks to Theorem \ref{thm:qpositive}, for any smooth strictly plurisubharmonic function $\psi$ on $U$ and any $\dbar$-closed $f\in \mathscr{D}_{(n,q)}(U)$, it is enough to show that there exists a solution $\dbar u=f$ satisfying 
\begin{equation}\label{eq:futoshiki}
\int_U |u|^2_{\omega_0}e^{-(\widetilde{\varphi}+\psi)}dV_{\omega_0}\leq \int_U \langle (c+[\ddbar \psi, \Lambda_{\omega_0}])^{-1}f, f \rangle_{\omega_0} e^{-(\widetilde{\varphi}+\psi)}dV_{\omega_0}.
\end{equation}
Consider the $\dbar$-closed $(n+m, q)$-form $\pi^\star f\wedge dW$ on $U\times D$.  
By assumption of $\varphi$, we can get a solution of $\dbar \widetilde{v}=\pi^\star f\wedge dW$ satisfying an $L^2$-estimates.
Take the $L^2$-minimal solution $\widetilde{u}$.  
``If" $\widetilde{u}$ has the form 
\begin{equation}\label{eq:minimal}
\widetilde{u}=\sum_{1\leq j_1<\cdots<j_{q-1}\leq n}\widetilde{u}_{j_1\cdots j_{q-1}}dZ\wedge dW\wedge d\bar{z}_{j_1}\wedge \cdots \wedge d\bar{z}_{j_{q-1}},
\end{equation}
each coefficient $\widetilde{u}_{j_1\cdots j_{q-1}}$ is holomorphic in $w$ and invariant under the rotation of $w$ due to the uniqueness of the minimal solution. Hence, we have $\widetilde{u}=\pi^\star u\wedge dW$ for some $u\in C^\infty_{(n, q-1)}(U)$ satisfying $\dbar u=f$ and the inequality (\ref{eq:futoshiki}) on $U$.  

However, the fact that $\widetilde{u}$ has the form (\ref{eq:minimal}) does not immediately follow from that $\widetilde{u}$ is the minimal $L^2$-solution, which is pointed out by Wang Xu. 
While there are still technical problems, we believe that Theorem \ref{thm:qpositive} is a valid way to solve Problem \ref{thm:qprekopa}. 
%In this section, we discuss the definition of uniform $q$-positivity for singular Hermitian metrics on line bundles. 

%Under this formulation, we can show Theorem \ref{thm:qprekopa} without assuming the condition that $D$ is bounded. 
%The proof remains the same. 

%\begin{theorem}$($cf. Theorem \ref{thm:qprekopa}$)$\label{thm:shmprekopa}
%	Let $U$ be a bounded domain in $\mathbb{C}^n_z$ and $D$ be a pseudoconvex domain in $\mathbb{C}^m_w$. 
	%Let $\varphi$ be a smooth function on $U\times D$. 
	%Assume that 
	%$\varphi$ satisfies the conditions $(1)$ and $(2)$ in Theorem \ref{thm:qprekopa}. We define the function $\widetilde{\varphi}$ on $U$ by 
	%$$
%	e^{-\widetilde{\varphi}(z)}:=\int_{w\in D}e^{-\varphi(z, w)}d\omega_1(w).
	%$$
	%Then, if $\widetilde{\varphi}\equiv -\infty$, $(\mathbb{C}, e^{-\widetilde{\varphi}})$ is uniformly $(q-1)$-$c$-positive in the sense of Definition \ref{def:shmpartial}. 
%\end{theorem}
As an application of Theorem \ref{thm:rc}, we also propose the following problem. 
The reason why Theorem \ref{thm:rc} is useful to Problem \ref{conj:directrc} is the same reason why Theorem \ref{thm:qpositive} is useful to Problem \ref{thm:qprekopa}.

\begin{problem}\label{conj:directrc}
	Assume that $E$ is weakly RC-positive. Is $S^kE\otimes \det E$ RC-positive for every $k\geq 1$?
\end{problem}

Weak RC-positivity of $E$ implies that $\mathscr{O}_E(1)$ is uniformly $(\dim X-1)$-positive, where $\mathscr{O}_E(1)$ is the tautological line bundle over the projectivized bundle $\mathbb{P}(E^\star)$ (cf. Proposition \ref{prop:uniform}). 
This problem asserts that if $\mathscr{O}_E(1)$ is uniformly $(\dim X-1)$-positive, $\pi_\star(K_{\mathbb{P}(E^\star)/X}\otimes \mathscr{O}(r+k))\cong S^kE\otimes \det E$ is RC-positive. 
This is related to the following conjecture raised by Yang. 

\begin{conjecture}$($\cite[Question 7.11]{Yan18}$)$\label{conj:yang}
	Assume that $E$ is weakly RC-positive.
	Then $E$ is RC-positive.
	\end{conjecture}

%%%%%%%%%%%%%%%%%%%%%%%%%%%%%%%%%%%%%%%%%%%%%%%%%%%%%%
%%%%%%%%%%%%%%%%%%%%%%%%%%%%%%%%%%%%%%%%%%%%%%%%%%%%%%

%%%%%%%%%%%% References %%%%%%%%%%%%%
%%
%<Author name> is written as Initial of Given Name, and Family Name.
%<Title> is written in roman letters.
%<Journal name> should be abbreviated according to
% the MR Serials Abbreviations List of Mathematical Reviews:
% (Abbreviations of Names of Serials; http://www.ams.org/mr-database)
%For <Pages>, use en-dash "--" between page numbers.
%%

\end{document}